# Formulas for $\pi(x)$ and the $n$th Prime

## Sebastian Martin Ruiz and Jonathan Sondow

**1. THE FORMULAS.** In [3], [4] the first author gave a formula for the $n$th prime number $p_n$ that involves only the elementary operations $+, -, \times, \div$ and the floor function $\lfloor \cdot \rfloor$. His conditional proof assumed certain inequalities based on the Prime Number Theorem. In this note, we use inequalities of Rosser and Schoenfeld [1] to give a complete proof of the slightly modified formula

$$(1) \qquad p_n = 2 + \sum_{k=2}^{\lfloor 2n \log n + 2 \rfloor} \left(1 - \left\lfloor \frac{\pi(k)}{n} \right\rfloor\right), \quad n > 1,$$

where $\pi(x)$, defined as the number of primes not exceeding $x$, is given by the formula

$$(2) \qquad \pi(x) = \sum_{j=2}^{\lfloor x \rfloor} \left(1 + \left\lfloor \frac{2 - \sum_{i=1}^{j}\left(\left\lfloor \frac{j}{i} \right\rfloor - \left\lfloor \frac{j-1}{i} \right\rfloor\right)}{j} \right\rfloor\right).$$

After the proof, we indicate various ways to modify and implement the formulas so that they operate in times $O((n \log n)^{3/2})$ and $O(x^{3/2})$, respectively.

**2. PROOF.** For $n$ a positive integer, let

$$(3) \qquad d(n) := \sum_{d \mid n} 1$$

denote the number of divisors of $n$. In [2], the first author found the formula

$$(4) \qquad d(n) = \sum_{i=1}^{n} \left( \left\lfloor \frac{n}{i} \right\rfloor - \left\lfloor \frac{n-1}{i} \right\rfloor \right),$$

which holds since the quantity in parentheses is 1 or 0 according as $i$ does or does not divide $n$.

Let $F$ be the characteristic function of the set of prime numbers

$$(5) \qquad F(n) := \begin{cases} 1 & \text{if } n \text{ is prime,} \\ 0 & \text{otherwise.} \end{cases}$$

From (3), we have $d(n) = 2$ if $n$ is prime, and $d(n) > 2$ if $n$ is composite. Since $2 \le d(n) \le n$ for $n > 1$, we have the formula

$$(6) \qquad F(n) = 1 + \left\lfloor \frac{2 - d(n)}{n} \right\rfloor, \quad n > 1.$$

Using (5), we write the prime-counting function $\pi(x)$ as the sum

$$(7) \qquad \pi(x) = \sum_{j=2}^{\lfloor x \rfloor} F(j)$$

with the convention that any sum $\sum_{i=a}^{b}$ is zero if $a > b$. From (7), (6), (4), we obtain formula (2) for $\pi(x)$.

In order to derive formula (1) for $p_n$ from (2), we require the following lemma.

**Lemma 1.** *For $n > 1$, we have the inequalities*

(8) $$\pi(2n \log n + 2) < 2n,$$

(9) $$p_n < 2n \log n + 2.$$

*Proof.* Rosser and Schoenfeld [1] proved that

(10) $$p_n > n \log n, \quad \text{all } n,$$

(11) $$p_n < n \log n + n(\log \log n - 1/2), \quad n > 20.$$

From (10), we have $p_{2n} > 2n \log 2n$. Since $\pi(p_{2n}) = 2n$, it follows that $\pi(2n \log 2n) < 2n$, which implies (8) if $n > 1$.

To prove (9) for $n > 1$, we verify it numerically for $n = 2, 3, \ldots, 20$, and note that (11) implies (9) for $n > 20$. This proves the lemma. ●

For $n > 1$, Lemma 1 implies that

$$\left\lfloor \frac{\pi(k)}{n} \right\rfloor = \begin{cases} 0 & \text{if } 1 \leq k \leq p_n - 1, \\ 1 & \text{if } p_n \leq k < 2n \log n + 2. \end{cases}$$

The desired formula for the $n$th prime number follows immediately. This completes the proof of (1) and (2). ●

**3. IMPROVEMENTS.** As they stand, the formulas for $p_n$, $\pi(x)$ and $d(n)$ operate in times $O((n\log n)^3)$, $O(x^2)$ and $O(n)$, respectively. We can improve these bounds by modifying the formulas as follows. If $i$ divides $n$, so does $n/i$, so to compute $d(n)$ it suffices to consider only $i \leq \sqrt{n}$. This reduces the times for $d(n)$ and $\pi(x)$ to $O(n^{1/2})$ and $O(x^{3/2})$, respectively. For $p_n$, compute $d(j)$ (and thus $F(j)$) for $j < 2n\log n + 2$, then compute $\pi(k)$ recursively as $\pi(k) = \pi(k-1) + F(k)$. This reduces the time for $p_n$ to $O((n\log n)^{3/2})$.

We can also improve the computation times (but not the O(.) bounds) in the following two ways. First, instead of the floor of $n/i$, use the integer quotient of $n$ by $i$. Second (as P. Sebah [5] has pointed out), in formulas (2) and (7) for $\pi(x)$, after $j=2$ we only need to sum over odd numbers, after $j=3$ only over numbers relatively prime to 6, and similarly for other moduli 30, 210, 2310,..., $m$. This "sieving" cuts computation time by a factor of $m/\varphi(m)$, where $\varphi$ is Euler's totient function.

We thank P. Sebah for discussions on the bounds, C. Rivera for the square root optimization, and J. McCranie for the quotient acceleration.

REFERENCES
1. J. B. Rosser and L. Schoenfeld, Approximate formulas for some functions of prime numbers, *Ill. J. Math.* **6** (1962) 64-94.
2. S. M. Ruiz, A functional recurrence to obtain the prime numbers using the Smarandache Prime Function, *Smarandache Notions Journal* **11** (2000) 56.
3. S. M. Ruiz, The general term of the prime number sequence and the Smarandache Prime Function, *Smarandache Notions Journal* **11** (2000) 59.
4. S. M. Ruiz, *Applications of Smarandache Functions and Prime and Coprime Functions*, American Research Press, Rehoboth, 2002.
5. P. Sebah, private communications, October-November, 2002.

*Avda. de Regla, 43 Chipiona 11550 Spain*
*smruiz@telefonica.net*

*209 West 97th St., New York, NY 10025*
*jsondow@alumni.princeton.edu*